\documentclass[10pt]{amsart}

\usepackage[latin1]{inputenc}

\usepackage{amsmath,amssymb,amsthm}
\usepackage[mathscr]{eucal}

\newtheorem{thm}{Theorem}[section]

\newtheorem{lemma}[thm]{Lemma}
\theoremstyle{definition}
\newtheorem*{prf}{Proof}
\newtheorem{definition}[thm]{Definition}

\def\bz{\bar{z}}

\newcommand\PD[2]{\frac{\partial #1}{\partial #2}}

\DeclareMathOperator{\extp}{ext}
\DeclareMathOperator{\intp}{int}
\DeclareMathOperator{\dbar}{\bar\partial}
\DeclareMathOperator{\dd}{\partial}
\def\ed{\mathrm{d}}
\def\wedgerum{\bigwedge^{0,*}(\mathbb{C}^d)^*}
\def\hrum{L_2(\mathbb{C}^d;\ed\lambda)\otimes \wedgerum}
\def\hrumw{L_2(\mathbb{C}^d;e^{-2W}\ed\lambda)\otimes \wedgerum}
\def\Hrum{\mathscr{H}}
\def\Hrumw{\mathscr{H}_W}

\def\b{\mathfrak{b}}
\def\pauli{\mathfrak{P}}
\def\dirac{\mathfrak{D}}
\def\paulif{\mathfrak{p}}
\def\tpauli{\widetilde{\mathfrak{P}}}
\def\tpaulif{\tilde{\mathfrak{p}}}

\let\phi=\varphi
\let\epsilon=\varepsilon
\numberwithin{equation}{section}

\title[Zero modes for the magnetic Pauli operator]{Zero modes for the magnetic\\ Pauli operator in even-dimensional\\ Euclidean space}
\author{Mikael Persson}
\address{Department of Mathematical Sciences\\
           Chalmers University of Technology\\
           G\"oteborg University\\
           SE-412 96 G\"oteborg, Sweden\\
}


\begin{document}
\maketitle

\begin{abstract}
We study the ground state of the Pauli Hamiltonian with a magnetic field in $\mathbb{R}^{2d}$. We consider the case where a scalar potential $W$ is present and the magnetic field $B$ is given by $B=2i\dd\dbar W$. The main result is that there are no zero modes if the magnetic field decays faster than quadratically at infinity. If the magnetic field decays quadratically then zero modes may appear, and we give a lower bound for the number of them. The results in this paper partly correct a mistake in a paper from 1993.
\end{abstract}

\section{Introduction and main result}

The Pauli operator $\pauli$ in $\mathbb{R}^n$
describes a charged spin-$\frac12$ particle in a
magnetic field. Along with the Dirac operator, it lies
in the base of numerous models in quantum physics. The
problem about zero modes, the bound states with zero
energy, is one of many questions to be asked about the
spectral properties of these operators.

Zero modes were discovered in~\cite{ac} in dimension
$n=2$. Unlike the purely electric interaction, a
compactly supported magnetic field can generate zero
modes, as soon as the total flux of the field is
sufficiently large. Quantitatively, this is expressed
by the famous Aharonov-Casher formula. The
two-dimensional case is by now quite well studied; the
AC formula is extended to rather singular magnetic
field, moreover, if the total flux is infinite (and
the field has constant sign), there are infinitely
many zero modes.

On the other hand, in the three-dimensional case the
presence of zero modes is a rather exceptional
feature, and the conditions for them to appear are not
yet found, see the discussion in~\cite{MeRo} and
references therein.

Even less clear is the situation in the higher
dimensions. In~\cite{sh}, for \emph{even} $n$ some
sufficient conditions for the infiniteness of the
number of zero modes  were found, requiring, in
particular, that the field decays rather slowly (more
slowly than $r^{-2}$) at infinity. On the other hand,
in~\cite{og}, again for even $n$, the case where a
finite number of zero modes should appear was
considered. Under the assumption of a rather regular
behavior of the \emph{scalar potential} of the
magnetic field at infinity the number of zero modes
was calculated. In particular, for a field with
compact support or  decaying faster than quadratically
at infinity the formula in~\cite{og} implies the
absence of zero modes, thus making a difference with
the two-dimensional situation.

Unfortunately, it turned out that the reasoning in~\cite{og} 
contains an error. A miscalculation in an
important integral leads to an erroneous conclusion,
thus destroying the  final results. This is the reason
for us to return to the question on zero modes in the
even higher-dimensional  case. We try to revive the
results in~\cite{og} and succeed partially.

We use the  representation of the Pauli and Dirac
operators in the terms of multi-variable complex
analysis proposed in \cite{sh} and used further in
~\cite{og}. This approach puts a certain restriction
on the class of magnetic field considered, equivalent
to the existence of  a scalar potential. At the moment
it is unclear how to treat the general case.

Under the above condition,  the operators are
represented as acting on the complex forms, the action
expressed via the $\dbar$ operator. The mistake in
~\cite{og} occurs in calculating the $L_2$ norm of the
form one gets after applying the $\dbar$ operator.  We
present the detailed analysis of this miscalculation
in Section~\ref{sec:log}.

The strategy of our treatment of zero modes differs
from the one in~\cite{ac} and other previous papers
including~\cite{og}. Usually, when studying zero
modes, one shows first that they, after having been
multiplied by some known factor, are holomorphic
function in the whole space; then one easily counts
the number of such functions. This strategy fails in
our case, so we use another one, involving more
advanced machinery of complex and real analysis. The
main ingredient of the proofs is a combination of the
techniques of using the Bochner-Martinelli-Koppelman
kernel to solve a $\dbar$ equation and the use of a
weighted Hardy-Littlewood-Sobolev inequality to
estimate that solution.

As a result, we establish some of the properties
presented in~\cite{og}. We show that there are no zero
modes if the magnetic field decays faster than
quadratically at infinity (in particular, if it is
compactly supported). Another result is that zero
modes may exist if the magnetic field decays exactly
quadratically, and the formula in~\cite{og} gives a
lower bound for their number.

\subsection{The Pauli operator}\label{sec:pauli}

Let $x=(x^1,\ldots,x^{2d})$ denote the usual Euclidean coordinates in $\mathbb{R}^{2d}$. According to the Maxwell equations, a magnetic field $B$ in $\mathbb{R}^{2d}$ is a real closed two-form
\begin{equation}\label{eq:magn}
B(x)=\sum_{j<k} b_{j,k}(x) \ed x^j\wedge \ed x^k.
\end{equation}
Throughout this paper we assume that all the coefficient functions $b_{j,k}$ belong to $C^\infty(\mathbb{R}^{2d})$. The condition that the magnetic field $B$ is closed is given by
\begin{equation*}
0=\ed B=\sum_{j<k<l} \left(\frac{\partial b_{j,k}}{\partial x^l}-\frac{\partial b_{j,l}}{\partial x^k}+\frac{\partial b_{k,l}}{\partial x^j}\right) \ed x^j\wedge \ed x^k \wedge \ed x^l,
\end{equation*}
where $\ed$ is the usual exterior differential operator. Since $B$ is closed there exists a one-form
\begin{equation*}
a(x)=\sum_{j=1}^n a_j(x) \ed x^j
\end{equation*}
satisfying
\begin{equation*}
B=\ed a=\sum_{j<k} \left(\frac{\partial a_k}{\partial x^j}-\frac{\partial a_j}{\partial x^k}\right)\ed x^j\wedge \ed x^k
\end{equation*}

Any such one-form $a$ is called a magnetic one-form or magnetic vector potential. It is not unique. In fact, given one magnetic one-form, another one is obtained by adding $\ed f$ for some regular function $f$. The choice of magnetic one-form $a$ is usually refered to as the choice of gauge.

The analysis of the Pauli operator was successful in~\cite{sh} using complex analysis under a condition that the magnetic field is a complex $(1,1)$-type form. It is not clear what this condition means physically, but to be able to use the theorey of complex analysis in several variables, we will throughout use the same assumption. Thus, the coefficient functions in~\eqref{eq:magn} of the closed $2$-form $B$ must satisfy the $d(d-1)$ equations
\begin{equation}\label{eq:ettett}
\begin{cases}
b_{2j-1,2k-1} = b_{2j,2k},& \\
b_{2j-1,2k} = - b_{2j,2k-1}& 
\end{cases}\qquad\text{for}\ 
j+1\leq k \leq d,\quad 1\leq j\leq d-1.
\end{equation}

The spinless Schr\"{o}dinger operator $H$ in $\mathbb{R}^{2d}$ corresponding to the magnetic field $B$ is defined in $L_2(\mathbb{R}^{2d})$ as
\begin{equation*}
H=\sum_{j=1}^{2d} \left(-i \frac{\partial}{\partial x^j}-a_j\right)^2.
\end{equation*}

We are interested in spin-$\frac12$ particles (including the electron). Such systems are described by the Pauli operator $\pauli$, acting in $L_2(\mathbb{R}^{2d})\otimes \mathbb{C}^{2^d}$. Let $\{\gamma^j\}_{j=0}^{2d}$ be Hermitian $2^d\times 2^d$ matrices satisfying
\begin{equation}\label{eq:clifford}
\gamma^j\gamma^k+\gamma^k\gamma^j = 2 \delta^{jk}I_{2^d},
\end{equation}
where $I_{2^d}$ denotes the $2^d\times 2^d$ identity matrix. These matrices $\{\gamma^j\}_{j=0}^{2d}$ generate a Clifford algebra, and are usually called the Dirac matrices. The Pauli operator $\pauli$ is defined by
\begin{equation*}
\pauli=HI_{2^d}+\sum_{0<j<k} i b_{jk}(x)\gamma^j\gamma^k.
\end{equation*}
To be more precise, $\pauli$ is first defined on $C^\infty_0\otimes \mathbb{C}^{2^d}$, where it is essential self-adjoint (see~\cite{ch}). We denote the self-adjoint closure by $\pauli$. The Pauli operator $\pauli$ can also be written as $\pauli=\dirac^2$, where $\dirac$ is the self-adjoint Dirac operator
\begin{equation*}
\dirac=\sum_{j=1}^{2d} \left(-i \frac{\partial}{\partial x^j}-a_j\right)\gamma^j.
\end{equation*}
From this it follows that the Pauli operator is non-negative. 

\subsection{The main result}

\begin{thm}\label{thm:main}
Assume that the equations in~\eqref{eq:ettett} are satisfied, and that there exist constants $C>0$ and $\rho>2$ such that
\begin{equation*}
|B(x)|\leq \frac{C}{(1+|x|)^\rho}\quad \text{for all } x\in\mathbb{R}^{2d}.
\end{equation*}
Then
\begin{equation*}
\dim\ker \pauli = 0.
\end{equation*}
\end{thm}
We will prove this theorem in Section~\ref{sec:bounded}. The case $|B(x)|\sim 1/|x|^2$ as $|x|\to\infty$ is more complicated. In section~\ref{sec:log} we give an example of a magnetic field satisfying
\begin{equation*}
|B(x)|=\frac{\Phi (d-1)}{2 |x|^2},\quad \text{for large values of $|x|$},
\end{equation*}
such that $\dim\ker \pauli=0$ if $|\Phi|<d$ and $\dim\ker \pauli > 0$ otherwise. This result is somehow strange and suggests that the situation for magnetic fields with a quadratic decay is quite complicated and unstable.

\subsection{Complex analysis and Differential forms}\label{sec:def}
Let us now switch to the complex analysis viewpoint. We identify $x=(x^1,\ldots,x^{2d})$ in $\mathbb{R}^{2d}$ with $z=(z^1,\ldots,z^d)$ in $\mathbb{C}^d$, where $z^j=x^{2j-1}+ix^{2j}$. We define tangent and cotangent vectors by
\begin{equation*}
\begin{aligned}
\frac{\partial}{\partial z^j} &= \frac12\left(\frac{\partial}{\partial x^{2j-1}}-i\frac{\partial}{\partial x^{2j}}\right),\\
\frac{\partial}{\partial \bz^j} &= \frac12\left(\frac{\partial}{\partial x^{2j-1}}+i\frac{\partial}{\partial x^{2j}}\right),\\
\ed z^j & = \ed x^{2j-1}+i \ed x^{2j},\quad\text{and}\\
\ed \bz^j & = \ed x^{2j-1}-i \ed x^{2j}.
\end{aligned}
\end{equation*}

Written in complex terms, the magnetic field $B$ can be written as a sum of $(1,1)$, $(2,0)$, and $(0,2)$ type forms as

\begin{equation*}
B(z) = \sum_{j,k=1}^d \b_{j,k}(z) \ed z^j \wedge \ed \bz^k + \sum_{j,k=1}^d \b'_{j,k}(z) \ed z^j \wedge \ed z^k + \sum_{j,k=1}^d \b''_{j,k}(z) \ed \bz^j \wedge \ed \bz^k.
\end{equation*}

The equations in~\eqref{eq:ettett} state that $B$ is of type $(1,1)$ which means that all coefficient functions $\b'_{j,k}$ and  $\b''_{j,k}$ in the representation above vanish, so the magnetic field $B$ has the form
\begin{equation}\label{eq:cettett}
B(z)=\sum_{j,k=1}^d \b_{j,k}(z) \ed z^j \wedge \ed \bz^k
\end{equation}

To magnetic fields that are $(1,1)$-type forms there exist scalar potentials $W\in C^\infty(\mathbb{C}^d\to\mathbb{R})$ satisfying
\begin{equation}\label{eq:magnetett}
 B=2i\dd\dbar W,
\end{equation}
see~\cite{wells}. In~\cite{sh} it was shown that the Dirac and Pauli operators can be defined in terms of $W$ and operators acting on differential forms in a very nice way. For the sake of completeness we show how this is done.

Let $\bigwedge^{0,q}(\mathbb{C}^d)^*$ denote the space of $(0,q)$-type differential forms and let $\wedgerum=\bigoplus_{q=0}^d\bigwedge^{0,q}(\mathbb{C}^d)^*$. The Dirac operator $\dirac$ is realized as an operator in the Hilbert space $\Hrum:=\hrum$ in the way
\begin{equation}\label{eq:kompldirac}
 \dirac=2\left(\dbar_W+\dbar_W^*\right).
\end{equation}
Here
\begin{equation*}
\begin{aligned}
 \dbar_W & = \dbar+\extp(\dbar W) = \sum_{j=1}^d\extp(\ed\bz^j)\left(\PD{}{\bz^j}+\PD{W}{\bz^j}\right),\\
 \dbar_W^* & = \dbar^*+\intp(\dd W) = \sum_{j=1}^d\intp(\ed z^j)\left(-\PD{}{z^j}+\PD{W}{z^j}\right),
\end{aligned}
\end{equation*}
$\extp(\ed\bz^j)$ is the operator on $\wedgerum$ acting as
\begin{equation*}
 \extp(\ed\bz^j)\eta = \ed\bz^j\wedge\eta,\quad \text{for}\ \eta\in\wedgerum,
\end{equation*}
and $\intp(\ed z^j)$ is the adjoint operator of $\extp(\ed\bz^j)$ in $\Hrum$.

To see that \eqref{eq:kompldirac} is true we use the anti-commutation relations

\begin{equation*}
 \begin{aligned}
  \mbox{}[\extp(\ed\bz^j),\extp(\ed\bz^k)]_+&=0;\\
  \mbox{}[\intp(\ed z^j), \intp(\ed z^k)]_+&=0;\\
  \mbox{}[\extp(\ed \bz^j), \intp(\ed z^k)]_+&=\delta^{jk}.
 \end{aligned}
\end{equation*}
By defining
\begin{equation*}
 \begin{aligned}
  \gamma^{2j-1}&=i(\extp(\ed \bz^j)-\intp(\ed z^j));\\
  \gamma^{2j}&=-(\extp(\ed\bz^j)+\intp(\ed z^j))
 \end{aligned}
\end{equation*}
one can easily check that
\begin{equation*}
 [\gamma^j,\gamma^k]_+ = 2\delta^{jk}.
\end{equation*}
Hence $\{\gamma^j\}$ so defined satisfies the relation~\eqref{eq:clifford} of a Clifford algebra. Now it is easy to see that~\eqref{eq:kompldirac} holds:

\begin{equation*}
 \begin{aligned}
  2 \left(\dbar_W+\dbar^*_W\right) & = \sum_{j=1}^d \left(-i\gamma^{2j-1}-\gamma^{2j}\right)\left(\frac{\partial}{\partial\bz^j}+\frac{\partial W}{\partial\bz^j}\right)
 +\left(i\gamma^{2j-1}-\gamma^{2j}\right)\left(-\frac{\partial}{\partial z^j}+\frac{\partial W}{\partial z^j}\right)\\
  & = \sum_{j=1}^d i\gamma^{2j-1}\left(-\frac{\partial}{\partial x^{2j-1}}-i\frac{\partial W}{\partial x^{2j}}\right)	
+ \gamma^{2j}\left(-i\frac{\partial}{\partial x^{2j}}-\frac{\partial W}{\partial x^{2j-1}}\right)\\
  & = \sum_{j=1}^{2d} \gamma^j\left(-i\frac{\partial}{\partial x^j}-a_j(x)\right) = \dirac
 \end{aligned}
\end{equation*}
where $a_{2j-1}(x)=\frac{\partial W}{\partial x^{2j}}$ and $a_{2j}(x)=-\frac{\partial W}{\partial x^{2j-1}}$, so $a=i(\dbar-\dd)W$, which fits well with~\eqref{eq:magnetett}, since $B=\ed a = (\dd+\dbar)a=(\dd+\dbar)i(\dbar-\dd)W=2i\dd\dbar W$.

For a form $\alpha\in\Hrum$ to belong to the kernel of $\pauli$ it is necessary and sufficient that $\alpha$ belongs to the kernel of the quadratic form

\begin{equation*}
 \paulif[\alpha]=4\int_{\mathbb{C}^d}\left(\left|\dbar_W\alpha\right|^2+\left|\dbar_W^*\alpha\right|^2\right)\ed\lambda(z),\quad \alpha\in\Hrum.
\end{equation*}

Let $U:\Hrum\to \Hrumw:=\hrumw$ be the unitary operator $U:\alpha\mapsto e^{W}\alpha$. Then $\pauli$ and $\tpauli=U\pauli U^*$ are unitarily equivalent. The quadratic form $\tpaulif$ on $\Hrumw$ corresponding to $\tpauli$ is given by
\begin{equation}\label{eq:wform}
 \tpaulif[\alpha] = 4\int_{\mathbb{C}^d}\left(\left|\dbar\alpha\right|^2+\left|\dbar^*\alpha\right|^2\right)e^{-2W}\ed\lambda(z),\quad \alpha\in\Hrumw.
\end{equation}
Here $\dbar^*$ is the adjoint operator to $\dbar$ in $\Hrumw$.

\section{Proof of Theorem~\ref{thm:main}}\label{sec:bounded}
\begin{lemma}
Under the same conditions as in Theorem~\ref{thm:main} there exists a scalar potential $W\in L^\infty(\mathbb{C}^d\to\mathbb{R})$ such that $2i\dd\dbar W = B$.
\end{lemma}

We know from~\cite{wells} that solutions $W$ exist, the essential part of this Lemma is that there exist a bounded solution to~\eqref{eq:magnetett}. 

\begin{prf}
First assume that $\rho$ is not an integer. Two applications of Theorem~9' in~\cite{boma} gives the existence a scalar potential $W$ and a constant $C>0$ such that
\begin{equation*}
|W(z)|\leq C\frac{1}{(1+|z|)^{\rho-2}}
\end{equation*}
which means that $W$ is bounded since $\rho>2$. If $\rho$ is an integer then we can replace $\rho$ by $\rho-1/2$ and use the same argument again.\qed
\end{prf}

To prove Theorem~\ref{thm:main} it is clearly enough to prove the following theorem.

\begin{thm}\label{thm:boundedW}
Assume that $W\in L^\infty(\mathbb{C}^d\to\mathbb{R})$. Then
\begin{equation*}
\dim\ker\pauli = 0.
\end{equation*}
\end{thm}

Since $\pauli$ and $\tpauli$ are unitarily equivalent we will show instead that $\dim\ker\tpauli=0$. We need some Lemmata.

\begin{lemma}\label{lem:haliso}
 Let $\Omega:\mathbb{C}^d\to \mathbb{C}$ be a homogeneous function of degree zero, and let $\Omega$ be bounded on the unit sphere $|z|=1$. Define the operator $T$ as
\begin{equation*}
 (T f)(z) = \int_{\mathbb{C}^d} \frac{\Omega(z-\zeta)}{|z-\zeta|^{2d-1}} f(\zeta)d\lambda(\zeta).
\end{equation*}
Then $T$ is bounded as an operator from $L_2(\mathbb{C}^d)$ to $L_{2d/(d-1)}(\mathbb{C}^d)$.
\end{lemma}

\begin{prf}
This is a special case of the Hardy-Littlewood-Sobolev theorem, see Theorem~V.1 in~\cite{st}.\qed
\end{prf}
Let $K_q(\zeta,z)$ be the Bochner-Martinelli-Koppelman kernel
\begin{equation*}
 K_q(\zeta,z) = \frac{(d-1)!}{2^{q+1}\pi^d}\frac{1}{|\zeta-z|^{2d}} \sum_{\substack{j,J\\|L|=q+1}} \epsilon_{jJ}^L (\overline{\zeta}^j-\bz^j)(*\ed\zeta^L)\wedge \ed\bz^J.
\end{equation*}
Here $J$ is a multiindex of length $q$ and if $A$ and $B$ are ordered subsets of $\{1,2,\ldots,d\}$ then $\epsilon_B^A$ denotes the sign of the permutation which takes $A$ into $B$ if $|A|=|B|$ and zero if $|A|\neq|B|$. If $A\subset \{1,2,\ldots,d\}$ and $|A|=q$ then
\begin{equation*}
 *\ed\zeta_A = \frac{(-1)^{q(q-1)/2}}{2^{d-q}i^d}\ed\zeta^A \wedge \left(\bigwedge_{\nu\in A'}\ed\overline\zeta^\nu\wedge \ed\zeta^\nu\right),
\end{equation*}
where $A'$ is the complementary multiindex of $A$. We see that $K_q(\zeta,z)$ is of type $(d,d-q-1)$ in $\zeta$ and $(0,q)$ in $z$.

\begin{lemma}\label{lem:beta}
Let $\alpha\in\Hrumw$ be a $(0,q)$-type form, $1\leq q\leq d-1$, satisfying $\dbar\alpha=\dbar^*\alpha=0$. Then the $(0,q-1)$-type form
\begin{equation}
 \beta(z)=-\int_{\mathbb{C}^d} \alpha(\zeta)\wedge K_{q-1}(\zeta,z)
\end{equation}
satisfies $\dbar\beta=\alpha$ in the sense of distributions. Moreover, there exists a constant $C>0$ such that
\begin{equation}
 \int_{R<|z|<2R} \frac{|\beta(z)|^2}{|z|^2}\ed\lambda(z) \leq C \|\alpha\|_{\Hrumw}^2
\end{equation}
for all $R>0$, where the constant $C$ does not depend on $\alpha$ or $R$.
\end{lemma}

\begin{prf}
 Let $\eta_k$ be a family of cut-off functions, such that $\eta_k(\zeta)=1$ if $|\zeta|<k$, $\eta_k(\zeta)=0$ if $|\zeta|>k+1$ and $|\dbar\eta_k|\leq 2$ for all $k=1,2,\ldots$. Since the form $\alpha$ belongs to the kernel of the elliptic Pauli operator (with smooth coefficient functions), it must itself be smooth. The forms $\eta_k\alpha$ are smooth and have compact support and thus, according to the Bochner-Martinelli-Koppelman formula (see~\cite{ra}), they satisfy
\begin{equation*}
\eta_k\alpha(z)=-\int \dbar(\eta_k \alpha)\wedge K_q(\zeta,z) - \dbar_z\int(\eta_k\alpha)\wedge K_{q-1}(\zeta,z),\quad \text{for $k=1,2,\ldots$.}
\end{equation*}
Let $\Phi$ be a test form with support in $|z|<M$. Fix $\epsilon>0$. Then
\begin{align*}
|\langle \beta, \dbar^*\Phi\rangle - \langle \alpha,\Phi \rangle | & = \left|\langle \beta, \dbar^*\Phi\rangle +\left\langle \int_\zeta (\eta_k\alpha)\wedge K_{q-1}(\zeta,z),\dbar^*\Phi\right\rangle \right.\\
&\qquad \left.+\left\langle \int_\zeta \dbar(\eta_k\alpha)\wedge K_q(\zeta,z),\Phi \right\rangle +\langle \eta_k\alpha, \Phi\rangle- \langle \alpha,\Phi \rangle \right|\\
& \leq \left|\left\langle \int_\zeta (\eta_k-1)\alpha\wedge K_{q-1}(\zeta,z),\dbar^*\Phi\right\rangle\right|\\
& \qquad + \left| \left\langle \int_\zeta \dbar\eta_k\wedge \alpha\wedge K_q(\zeta,z),\Phi \right\rangle \right|\\
& \qquad + \left|\langle (\eta_k-1)\alpha,\Phi \rangle \right|\\
& = I_1+I_2+I_3.	
\end{align*}
We will let $k$ tend to infinity. For $k>2M$ we have $|K_{q-1}(\zeta,z)|\leq C |\zeta|^{1-2d}$. We get
\begin{align*}
I_1 &\leq \int_{|z|<M}\int_{|\zeta|>k} |\eta_k(\zeta)-1|\cdot |\alpha(\zeta)|\cdot |K_{q-1}(\zeta,z)|\ed\lambda(\zeta)\ |\dbar^*\Phi(z)|\,\ed\lambda(z)\\
 & \leq C \sup |\dbar^*\Phi|\cdot \|\alpha\|_{\Hrum} \int_{|z|<M}  \left(\int_{|\zeta|>k} |K_{q-1}(\zeta,z)|^2 \ed\lambda(\zeta) \right)^{1/2}\ed\lambda(z)\\
&\leq C \sup |\dbar^*\Phi|\cdot\|\alpha\|_{\Hrum} \left(\int_k^\infty r^{2-4d+2d-1}\ed r\right)^{1/2}\\
&\leq C \sup |\dbar^*\Phi|\cdot\|\alpha\|_{\Hrumw}\cdot k^{1-d}
\end{align*}
so $I_1 < \epsilon$ if $k$ is large enough. Similarily, for $I_2$, we have
\begin{align*}
I_2 &\leq \int_{|z|<M}\int_{k<|\zeta|<k+1} |\dbar \eta_k|\cdot |\alpha|\cdot |K_q(\zeta,z)|\ed\lambda(\zeta) \ |\Phi(z)|\,\ed\lambda(z)\\
& \leq C \sup |\Phi|\cdot \|\alpha\|_{\Hrum} \int_{|z|<M} \left( \int_{k<|\zeta|<k+1} |K_q(\zeta,z)|^2 e^{2W(\zeta)} \ed\lambda(\zeta)\right)^{1/2} \ed\lambda(z)\\
& \leq C \sup |\Phi|\cdot \|\alpha\|_{\Hrum} \left(\int_k^{k+1} r^{2-4d+2d-1}\ed r\right)^{1/2}\\
& \leq C \sup |\Phi|\cdot \|\alpha\|_{\Hrumw}\cdot k^{\frac12-d},
\end{align*}
so $I_2<\epsilon$ if $k$ is large enough. $I_3$ is equal to zero if $k$ is large enough, since then $(1-\eta_k)$ and $\Phi$ has disjoint support. We conclude that $\dbar\beta=\alpha$ in the sense of distributions.

To show the estimate~\eqref{eq:hjalpgen}, we use Lemma~\ref{lem:haliso}. Indeed, note that $\beta$ can be written as
\begin{equation*}
 \beta=\sum_{J} T_J\alpha_J
\end{equation*}
where $\alpha=\sum_J \alpha_J d\bz^J$, $|J|=q$, and all operators $T_J$ are of the kind in Lemma~\ref{lem:haliso}. Denote by $E_R$ the set $\{z\in\mathbb{C}^d\ :\ R<|z|\leq 2R\}$. Using the H\"older inequality and Lemma~\ref{lem:haliso}, we have
\begin{align*}
\int_{E_R} \frac{|\beta|^2}{|z|^2}\ed\lambda(z) &\leq \left(\int_{E_R} \frac{1}{|z|^{2d}}\ed\lambda(z)\right)^{1/d} \left(\int_{E_R} |\beta|^{2d/(d-1)}\ed\lambda(z)\right)^{(d-1)/d}\\
& \leq	 C \|\alpha\|_{\Hrum}^2 \leq	 C \|\alpha\|_{\Hrumw}^2.
\end{align*}
Note that the integral
\begin{equation*}
\int_{R<|z|<2R} \frac{1}{|z|^{2d}}d\lambda(z) = c_d \int_{R}^{2R} r^{-2d+2d-1}dr =c_d\log(2)
\end{equation*}
is independent of $R$, so the constant $C$ above is also independent of $R$.\qed
\end{prf}

\begin{prf}[of Theorem~\ref{thm:boundedW}]
Let $1\leq q \leq d-1$. Assume that $\alpha\in\Hrumw$ is a $(0,q)$-type form in the kernel of $\tpauli$. Then $\dbar\alpha=\dbar^*\alpha=0$, and so we get the form $\beta$ from Lemma~\ref{lem:beta}. We don't know, apriori, that $\beta$ belongs to the domain of the $\dbar$ operator. We introduce a family of cut-off functions to be able to integrate by parts.

Let $\phi_k(r)$, $k=1,2,\ldots$, be a $C^\infty$ family of cut-off functions, such that $\phi_k(r)=1$ if $0<r\leq 2^{k}$, $\phi_k(r)=0$ if $r\geq 2^{k+1}$ and such that $0\leq \phi_k$ and $|\phi_k'(r)|\leq 2^{1-k}$. Let $\chi_k(z)=\phi_k(|z|)$. We have
\begin{align*}
 0&=\langle \dbar^*\alpha,\chi_k\beta\rangle_{\Hrumw}\\
 &= \langle \alpha,\dbar(\chi_k\beta)\rangle_{\Hrumw}\\
 &= \int |\alpha|^2\chi_k e^{-2W}\ed\lambda + \int \alpha \cdot \overline{\dbar\chi_k\wedge \beta} e^{-2W}\ed\lambda\\
 &= I_k+II_k.
\end{align*}
The integration by parts above is permitted thanks to the cut-off function $\chi_k$. It is clear that $I_k\to ||\alpha||^2_{\Hrumw}$ as $k\to\infty$. We shall prove that $II_k\to 0$ as $k\to\infty$.

Let $m_k^2=\int_{E_k} |\alpha|^2 e^{-2W}d\lambda$. Then it holds that $\sum_{k}m_k^2=||\alpha||^2_{\Hrumw}<\infty$ so $m_k\to 0$ as $k\to\infty$. Since $\dbar\chi_k$ has support in $E_k$ and $|\dbar \chi_k|\leq C2^{-k}$ we have
\begin{align*}
|II_k| &\leq \int_{E_k} |\alpha|\cdot |\beta|\cdot |\dbar\chi_k|e^{-2W}\ed\lambda \\
&\leq C 2^{-k} \left(\int_{E_k} |\alpha|^2 \ed\lambda \right)^{1/2}\left(\int_{E_k} |\beta|^2 \ed\lambda\right)^{1/2}\\
& \leq C m_k\left(\int_{E_k} \frac{|\beta|^2}{|z|^2} \ed\lambda \right)^{1/2}\\
& \leq C  m_k \|\alpha\|_{\Hrumw}\to 0,\quad\text{as $k\to\infty$.}
\end{align*}
The first inequality is just the triangle inequality. The second one is the inequality for $\chi_k$ and the Cauchy-Schwarz inequality. In the third inequality we use the fact that $|z|\approx 2^{k}$, and in the fourth we use Lemma~\ref{lem:beta}.

Next let $q=0$, and assume that $\alpha$ is a $(0,0)$-type form in the kernel of $\tpauli$. According to~\eqref{eq:wform} $\alpha$ has to be an entire function in $z^1,\ldots,z^d$. Since the function $\alpha$ also belongs to $L_2(\mathbb{C}^d,e^{-2W}\ed\lambda)$ a Liouville-type argument gives that it must be zero.

Finally let $q=d$. Then~\eqref{eq:wform} implies that $\dbar^*\alpha=0$. If $\alpha=\hat\alpha \ed\bz^1\wedge \cdots \wedge \ed\bz^d$, then this means that
\begin{equation*}
 -\frac{\partial\hat\alpha}{\partial z_j}+2\frac{\partial W}{\partial z_j}\hat\alpha = 0,\qquad j=1,\ldots,d.
\end{equation*}
If we put $f(z)=e^{-2W(z)}\hat\alpha(z)$ we obtain
 \begin{equation*}
 \frac{\partial f}{\partial z_j} = 0,\qquad j=1,\ldots,d,
\end{equation*}
that is the function $f$ is an entire function in $\bz^1,\ldots,\bz^d$. Moreover the function $f$ belongs to $L_2(\mathbb{C}^d,e^{2W}\ed\lambda)$ so it must be zero.
\qed
\end{prf}

\section{Quadratically decaying magnetic fields}\label{sec:log}
The case of determining the kernel of the Pauli operator for potentials with a logarithmic growth, which includes quadratically decaying magnetic fields, is more complicated. Given a real number $\Phi$, denote by $N_d(\Phi)$ the number of all monomials in $d$ variables with degree less than $|\Phi|-d$. The following Theorem was proposed in~\cite{og}.

\begin{thm}\label{thm:mainthm}
 Assume that $W\in C^\infty(\mathbb{C}^d\to\mathbb{R})$ and that there exists a real constant $\Phi$ such that the limit
 \begin{equation*}
  \lim_{|z|\to\infty} \frac{e^{-W(z)}}{|z|^{\Phi}}
 \end{equation*}
exists and is greater than zero. Then
\begin{equation*}
 \dim\ker\pauli = N_d(\Phi).
\end{equation*}
\end{thm}

Let us sketch the idea of the proof in the case $d=2$. First, assume that $\Phi>0$, and that
\begin{equation*}
 \alpha = \alpha_{00}+\alpha_{10}\ed\bz^1+\alpha_{01}\ed\bz^2+\alpha_{11}\ed\bz^1 \wedge \ed\bz^2
\end{equation*}
is an element of $\ker\tpauli$. Then 

\begin{equation*}
\dbar\alpha = \frac{\partial\alpha_{00}}{\partial \bz^1} \ed\bz^1 + \frac{\partial\alpha_{00}}{\partial \bz^2} \ed\bz^2 + \left(\frac{\partial\alpha_{01}}{\partial \bz^1} - \frac{\partial\alpha_{10}}{\partial \bz^2}\right)\ed\bz^1 \wedge \ed\bz^2 
\end{equation*}
and thus

\begin{align*}
0 &=\int_{\mathbb{C}^2}|\bar\partial \alpha|^2 e^{-2W} \ed\lambda\\
  &= \int_{\mathbb{C}^2}\left(\left|\frac{\partial\alpha_{00}}{\partial \bz^1}\right|^2+\left|\frac{\partial\alpha_{00}}{\partial \bz^2}\right|^2+\left|\frac{\partial\alpha_{01}}{\partial \bz^1}-\frac{\partial\alpha_{10}}{\partial \bz^2}\right|^2\right)e^{-2W} \ed\lambda.
\end{align*}
However, in~\cite{og} this is written as

\begin{equation}\label{eq:wrong}
 \begin{aligned}
0 &=\int_{\mathbb{C}^2}|\bar\partial \alpha|^2 e^{-2W} \ed\lambda\\
  &= \int_{\mathbb{C}^2}\left(\left|\frac{\partial\alpha_{00}}{\partial \bz^1}\right|^2+\left|\frac{\partial\alpha_{00}}{\partial \bz^2}\right|^2+\left|\frac{\partial\alpha_{01}}{\partial \bz^1}\right|^2+\left|\frac{\partial\alpha_{10}}{\partial \bz^2}\right|^2\right)e^{-2W} \ed\lambda,
\end{aligned}
\end{equation}
which is not correct. The rest of the proof uses~\eqref{eq:wrong} and some arguments to show that $\alpha_{00}$, $\alpha_{10}$ and $\alpha_{01}$ must vanish. Then it is shown, correctly, that the term $\alpha_{11}\ed\bz^1\wedge\ed\bz^2$ contains elements in the kernel if $\Phi$ is big enough. It is similar if $\Phi<0$.

So, we know from~\cite{og} that if $W\sim -\Phi\log|z|$, as $|z|\to\infty$, for $|\Phi|>d$, then the kernel is non-empty, and the dimension of the kernel is at least $N_d(\Phi)$. We are not able to prove the Theorem proposed in~\cite{og}, but we can show the following Theorem.
\begin{thm}\label{thm:mainlog}
 Assume that the limit
 \begin{equation*}
  \lim_{|z|\to\infty}\frac{e^{-W(z)}}{|z|^\Phi}
 \end{equation*}
exists and is positive. If $|\Phi|<d$ then $\dim\ker\pauli=0$. If $|\Phi|\geq d$ then $\dim\ker\pauli\geq N_d(\Phi).$
\end{thm}

The proof goes in the same way as the proof of Theorem~\ref{thm:boundedW} so we will just point out the main differences. First of all we can assume that $\Phi\geq 0$. If $\Phi$ is negative we can apply a unitary transform that changes the sign of $W$.

We need a replacement of Lemma~\ref{lem:beta} where weights of polynomial growth are allowed. To prepare for this we introduce the Muckenhoupt weight class.
\begin{definition}
 A non-negative function $\psi$ is said to belong to the Muckenhoupt class $A(p,q)$, $1<p,q<\infty$, if there exists a constant $C>0$ such that
 \begin{equation*}
  \sup_{B\subset \mathbb{C}^d}\left(\frac{1}{|B|}\int_B \psi^q d\lambda(z)\right)^{1/q} \left(\frac{1}{|B|}\int_B \psi^{-p/(p-1)} d\lambda(z)\right)^{(p-1)/p}\leq C.
 \end{equation*}
Here the supremum is taken over all balls in $\mathbb{C}^d$ and $|B|$ denotes the Lebesgue measure of the ball $B$.
\end{definition}

\begin{lemma}\label{lem:vikten}
If $W$ satisfies the properties in Theorem~\ref{thm:mainlog} then the weight function $e^{-W}$ belongs to the Muckenhoupt class $A(2,2d/(d-1))$.
\end{lemma}

\begin{prf}
Let $\gamma=2d/(d-1)$. We should show that
 \begin{equation*}
  I:=\left(\frac{1}{|B|}\int_B e^{-\gamma W} d\lambda(z)\right)^{1/\gamma} \left(\frac{1}{|B|}\int_B e^{2 W} d\lambda(z)\right)^{1/2}\leq C,
 \end{equation*}
where $C$ does not depend on the ball $B$. From the assumptions on $e^{-W}$ we know that there exist positive constants $c_1$, $c_2$, $c_3$ and $c_4$ such that 
\begin{equation}\label{eq:est1}
 c_1|z|^\Phi \leq e^{-W(z)} \leq c_2|z|^\Phi,\qquad \text{if } |z|\geq 1
\end{equation}
and
\begin{equation}\label{eq:est2}
 c_3 \leq e^{-W(z)} \leq c_4,\qquad \text{if } |z|< 5.
\end{equation}
We divide the balls into different classes. Say that a ball $B=B(z_0,R)$ is of Type 1 if $|z_0|>3/2 R$ and otherwise of Type 2.

First, assume that $B$ is of Type 1. Then for $z\in B$ we have $|z|\leq |z_0|+R\leq 5/3 |z_0|$ and $|z|\geq |z_0|-R\geq 1/3 |z_0|$. If $|z_0|\geq 3$ we can use \eqref{eq:est1} and get
\begin{align*}
 I &\leq C \left(\frac{1}{|B|}\int_B |z|^{\gamma\Phi} d\lambda(z)\right)^{1/\gamma} \left(\frac{1}{|B|}\int_B \frac{1}{|z|^{2\Phi}}\ d\lambda(z)\right)^{1/2}\\
 & \leq C\left(\frac{1}{|B|}\int_B |z_0|^{\gamma\Phi} d\lambda(z)\right)^{1/\gamma} \left(\frac{1}{|B|}\int_B \frac{1}{|z_0|^{2\Phi}} d\lambda(z)\right)^{1/2}\\
 & = C \left(|z_0|^{\gamma\Phi}\right)^{1/\gamma} \left(\frac{1}{|z_0|^{2\Phi}}\right)^{1/2}\\
 & = C.
\end{align*}
If $|z_0|\leq 3$ then $|z|\leq 5$, so we can easily use~\eqref{eq:est2} to get that $I\leq C$ independent of $R$.

Now assume that $B$ is of Type 2. Then $B\subset B':=B(0,3R)$. Since $B$ and $B'$ are of comparable size, we have
\begin{equation}
 I\leq C\left(\frac{1}{|B'|}\int_{B'} e^{-\gamma W} d\lambda(z)\right)^{1/\gamma} \left(\frac{1}{|B'|}\int_{B'} e^{2 W} d\lambda(z)\right)^{1/2} =: J.
\end{equation}
If $R\leq 5/3$ we can use~\eqref{eq:est2} to get that $J\leq C$ independent of $R$. If $R>5/3$ we have
\begin{align*}
 J & \leq C\left(\frac{1}{R^{2d}}\left(\int_{|z|<5} (1/c_4)^\gamma d\lambda(z)+\int_{5<|z|<3R} \frac{|z|^{\gamma\Phi}}{c_2^\gamma}d\lambda(z)\right)\right)^{1/\gamma} \times\\
 & \qquad \left(\frac{1}{R^{2d}}\left(\int_{|z|<5} c_3^\gamma d\lambda(z)+\int_{5<|z|<3R} \frac{c_1^2}{|z|^{2\Phi}}d\lambda(z)\right)\right)^{1/2}
\end{align*}
In this product the first factor is $O(R^{\Phi})$ while the second factor is $O(R^{-\min(d,\Phi)})$ as $R\to\infty$. Since the expression clearly is bounded for bounded values of $R$ there exists a constant $C$ such that $J\leq C$ independent of $R$.

We conclude that $e^{-W}\in A(2,2d/(d-1))$.\qed
\end{prf}

The following Lemma replaces Lemma~\ref{lem:haliso}.
\begin{lemma}\label{lem:dinglu}
 Let $\Omega:\mathbb{C}^d\to \mathbb{C}$ be a homogeneous function of degree zero, and let $\Omega$ be bounded on the unit sphere $|z|=1$. Define the operator $T$ as
\begin{equation*}
 (T f)(z) = \int_{\mathbb{C}^d} \frac{\Omega(z-\zeta)}{|z-\zeta|^{2d-1}} f(\zeta)d\lambda(\zeta).
\end{equation*}
If the weight function $\psi$ belongs to the Muckenhoupt class $A\left(2,2d/(d-1)\right)$ then there exists a constant $C>0$, independent of $f$, such that
\begin{equation*}
\left(\int_{\mathbb{C}^d} \left((Tf)(z)\cdot \psi(z)\right)^{2d/(d-1)}\right)^{(d-1)/(2d)} \leq C \left(\int_{\mathbb{C}^d}\left|f(z)\psi(z)\right|^2\right)^{1/2}.
\end{equation*}
\end{lemma}

\begin{prf}
This is a special case of Theorem~1 in \cite{dilu}.\qed
\end{prf}

Finally we get the result that replaces Lemma~\ref{lem:beta}.

\begin{lemma}\label{lem:boundgen}
Let $\alpha\in\Hrumw$ be a $(0,q)$-type form, $1\leq q\leq d-1$, satisfying $\dbar\alpha=\dbar^*\alpha=0$. Then the $(0,q-1)$-type form
\begin{equation}\label{eq:beta}
 \beta(z)=-\int_{\mathbb{C}^d} \alpha(\zeta)\wedge K_{q-1}(\zeta,z)
\end{equation}
satisfies $\dbar\beta=\alpha$ in the sense of distributions. Moreover, there exists a constant $C>0$ such that
\begin{equation}\label{eq:hjalpgen}
 \int_{R<|z|<2R} \frac{|\beta(z)|^2}{|z|^2}e^{-2W(z)}d\lambda(z) \leq C \|\alpha\|_{\Hrumw}^2
\end{equation}
for all $R>0$, where the constant $C$ does not depend on $\alpha$ or $R$.
\end{lemma}

\begin{prf}
The part that $\beta$ solves $\dbar\beta=\alpha$ is just the same as in the proof of Lemma~\ref{lem:beta}. For the estimate we can use Lemma~\ref{lem:dinglu} to get
\begin{equation*}
\begin{aligned}
 \int_{E_R} |\beta|^2e^{-2W}\ed\lambda(z) &\leq \left(\int_{E_R}|z|^{-2d}\right)^{1/d}\left(\int_{E_R}\left(|\beta|e^{-W}\right)^{2d/(d-1)}\ed\lambda(z)\right)^{(d-1)/d}\\
 &\leq C \|\alpha\|^2_{\Hrumw}
 \end{aligned}
\end{equation*}
\qed
\end{prf}

\begin{prf}[of Theorem~\ref{thm:mainlog}]
First, let $1\leq q\leq d-1$. The proof runs in the same way as in the proof of Theorem~\ref{thm:boundedW}, but with the use of Lemmata with weights.

Next, let $q=0$, and assume that $\alpha$ is a $(0,0)$-type form in the kernel of $\tpauli$. According to~\eqref{eq:wform} $\alpha$ has to be an entire function in $z_1,\ldots,z_d$. Also belonging to $L_2(\mathbb{C}^d,e^{-2W}d\lambda)$, it must tend to zero at infinity. Hence it must be constant equal to zero by a Liouville type argument.

Finally, let $q=d$. Then~\eqref{eq:wform} implies that $\dbar^*\alpha=0$. If $\alpha=\hat\alpha d\bz_1\wedge \cdots \wedge d\bz_d$, then this means that the function $f(z)=e^{-2W(z)}\hat\alpha(z)$ is an entire function in $\bz_1,\ldots,\bz_d$. Moreover there exist constants $c_1$ and $c_2$ such that
 \begin{equation*}
\frac{c_1}{|z|^\Phi}\leq e^{W(z)} \leq \frac{c_2}{|z|^\Phi}
\end{equation*}
if $|z|$ is large enough. Now the condition that $\alpha\in\Hrumw$ means that $e^W f\in L_2(\mathbb{C}^d)$. This is the case if and only if $f$ is a polynomial in $\bz_1,\ldots,\bz_d$ of degree strictly less than $\Phi-d$. The dimension of the space of such polynomials is exactly $N_d(\Phi)$.
\qed
\end{prf}

\def\cprime{$'$}

\end{document}